\documentclass[11pt]{article}

\newcommand{\ssage}{{\em Sage\ }}

\usepackage{amsmath}
\usepackage{url}
\usepackage{longtable}
\usepackage[margin=1in]{geometry}
\usepackage{array}
\newcommand{\sageprog}{{\em Sage}}

\usepackage{graphicx}
\usepackage{pdflscape}

\usepackage{sagetex}

\title{Table of minimum ranks of graphs of order at most 7 and
  selected optimal matrices}

\author{Laura DeLoss\thanks{Department of Mathematics, Iowa State
    University, Ames, IA 50011, USA (delolau@iastate.edu,
    grout@iastate.edu, tmckay16@iastate.edu, smithj@iastate.edu,
    gtims@iastate.edu).} \and Jason Grout\footnotemark[1] \and Leslie
  Hogben\thanks{Department of Mathematics, Iowa State University,
    Ames, IA 50011, USA (lhogben@iastate.edu) and American Institute
    of Mathematics, 360 Portage Ave, Palo Alto, CA 94306, USA
    (hogben@aimath.org).} \and Tracy McKay\footnotemark[1] \and Jason
  Smith\footnotemark[1]\and Geoff Tims\footnotemark[1]}

\begin{document}
\maketitle

\begin{abstract}
  The minimum rank of a simple graph $G$ is defined to be the smallest
  possible rank over all symmetric real matrices whose $ij$th entry
  (for $i\neq j$) is nonzero whenever $\{i,j\}$ is an edge in $G$ and
  is zero otherwise.  Minimum rank is a difficult parameter to
  compute.  However, there are now a number of known reduction
  techniques and bounds that can be programmed on a computer; we have
  developed a program using the open-source mathematics software
  \ssage~to implement several techniques.  We have also established
  several additional strategies for computation of minimum rank.
  These techniques have been used to determine the minimum ranks of
  all graphs of order 7.  This paper contains a list of minimum ranks
  for all graphs of order at most 7.  We also present selected optimal
  matrices.
\end{abstract}


\section{Table of minimum ranks}

In this report, we give a table of minimum ranks and various bounds on
the minimum ranks for all graphs with at most 7 vertices.  For
selected graphs, we also list optimal matrices that realize the
minimum rank.  With only a few exceptions (which are italicized), the
information in this table was computed using the program \cite{SAGE},
written in the \sageprog\ system \cite{sage}.

For a computer-readable version of the data in
Table~\ref{tab:minrank}, see the \verb|data/| directory accompanying
this paper.  The \verb|data/minrank-without-precomputed.csv|
file is a comma-separated file that can be easily opened with a
spreadsheet program.

In Table~\ref{tab:minrank}, we have abbreviated the column headings
for convenient presentation.  An ``LB'' in the column heading
signifies that the column contains a lower bound for the minimum rank,
while an ``UB'' signifies an upper bound.  The explanations of the
column headings are:
\begin{description}
\item[Atlas \#] The Atlas of Graphs \cite{atlas} number of a graph $G$
\item[$|G|$] The number of vertices in the graph
\item[Size] The number of edges in the graph
\item[$\text{mr}(G)$] The minimum rank of the graph.  If the minimum rank was found by hand, then the number is italicized.
\item[LB] The maximum lower bound found by the program
\item[UB] The minimum upper bound found by the program
\item[Con.] ``T'' if the graph is connected, ``F'' otherwise.  Most of
  the rest of the bounds are only listed if the graph is connected.
\item[ZFS LB] The lower bound obtained from the minimal zero forcing
  set size
\item[Diam LB] The lower bound obtained from the diameter of the graph
\item[CC UB] The upper bound obtained from the clique cover number of
  the graph
\item[NP UB] The upper bound obtained from the fact that the graph is
  not planar (only listed for nonplanar graphs)
\item[NOP UB] The upper bound obtained from the fact that the graph is
  not outer planar (only listed for graphs that are not outer planar)
\item[Path UB] The upper bound obtained from the fact that the graph
  is not a path (only listed for graphs that are not paths)
\item[IS] ``T'' if the graph contains an induced forbidden subgraph
  for minimum rank 2, ``F'' otherwise.  If ``T'', then the graph has
  minimum rank at least 3; if ``F'', then the graph has minimum rank
  at most 2.
\item[CV] ``T'' if the graph has a cut vertex, ``F'' otherwise
\item[Tree] ``T'' if the graph is a tree, ``F'' otherwise
\end{description}

\newlength{\zfscol}\settowidth{\zfscol}{ZFS LB}
\newlength{\diamcol}\settowidth{\diamcol}{Diam LB}
\newlength{\cccol}\settowidth{\cccol}{CC UB}
\newlength{\npcolumn}\settowidth{\npcolumn}{NP UB}
\newlength{\nopcol}\settowidth{\nopcol}{NOP UB}
\newlength{\npcol}\settowidth{\npcol}{Path UB}

\begin{landscape}
\begin{longtable}{%
|c
|c
|c
|c
|c
|c
|c
|c
|c
|c
|c
|c
|c
|c
|c
|c|
}
  \caption{Table of minimum ranks}\label{tab:minrank}\\
\hline
Atlas \# & $|G|$ & Size & $\mathrm{mr}(G)$ & LB & UB & Con. & \parbox[b]{\zfscol}{\strut\centering ZFS LB} & \parbox[b]{\diamcol}{\strut\centering Diam LB} & \parbox[b]{\cccol}{\centering CC UB} & \parbox[b]{\npcolumn}{\centering NP UB} & \parbox[b]{\nopcol}{\centering NOP UB} & \parbox[b]{\npcol}{\centering Path UB} & IS & CV & Tree \\\hline\hline\endfirsthead
\hline
Atlas \# & $|G|$ & Size & $\mathrm{mr}(G)$ & LB & UB & Con. & \parbox[b]{\zfscol}{\strut\centering ZFS LB} & \parbox[b]{\diamcol}{\strut\centering Diam LB} & \parbox[b]{\cccol}{\centering CC UB} & \parbox[b]{\npcolumn}{\centering NP UB} & \parbox[b]{\nopcol}{\centering NOP UB} & \parbox[b]{\npcol}{\centering Path UB} & IS & CV & Tree \\\hline\hline\endhead
\hline\endlastfoot
  1 & 1 & 0 & 0 & 0 & 0 & T & 0 & 0 & 0 &  &  &  & F & F & T \\
2 & 2 & 0 & 0 & 0 & 0 & F &  &  &  &  &  &  &  & F & F \\
3 & 2 & 1 & 1 & 1 & 1 & T & 1 & 1 & 1 &  &  &  & F & F & T \\
\hline
4 & 3 & 0 & 0 & 0 & 0 & F &  &  &  &  &  &  &  & F & F \\
5 & 3 & 1 & 1 & 1 & 1 & F &  &  &  &  &  &  &  & F & F \\
6 & 3 & 2 & 2 & 2 & 2 & T & 2 & 2 & 2 &  &  &  & F & T & T \\
\hline
7 & 3 & 3 & 1 & 1 & 1 & T & 1 & 1 & 1 &  &  & 1 & F & F & F \\
8 & 4 & 0 & 0 & 0 & 0 & F &  &  &  &  &  &  &  & F & F \\
\end{longtable}
\end{landscape}

\section{Selected optimal matrices}

There are some graphs in Table~\ref{tab:minrank} for which the program
was not able to determine the minimum rank.  For most of these graphs,
we now list in Table~\ref{tab:optimal} optimal matrices which show
that the lower bound found by the program is indeed the minimum rank
(i.e., the optimal matrices listed in Table~\ref{tab:optimal} are
evidence that the italicized numbers in Table~\ref{tab:minrank} are
correct).  The minimum ranks of the remainder of the graphs (graphs
558, 669, 678, 679, 791, 1086, and 1135) are determined in
\cite[Proposition~4.1]{SMALLGRAPHS}.

For a computer-readable version of the data in
Table~\ref{tab:optimal}, see the \verb|data/| directory accompanying
this paper.  The \verb|data/minrank-witnesses.sage| is the \ssage source
code that contains the optimal matrices and can be loaded into a
running \ssage session or pasted into a \ssage notebook cell, while the
\verb|data/Matrix_Witnesses.sws| file is an equivalent \ssage worksheet
that can be uploaded and used in the \ssage notebook.

\begin{sagesilent}
import networkx.generators.atlas
atlas_graphs = [Graph(i) for i in networkx.generators.atlas.graph_atlas_g()]
def check_witness(Matrix,Number):
    print "Rank is ", Matrix.rank()
    print "Symmetric?", Matrix.is_symmetric()
    print "Isomorphic to graph G",Number," ",Graph(Matrix).is_isomorphic(atlas_graphs[Number])
    print Matrix
    show(Graph(Matrix))
\end{sagesilent}

\begin{sagesilent}
G={}

G[721]=Matrix([[-1, 1, 0, 0, 1,0,1],[1,-1,1,0,0,1,0],[0,1,-1,1,0,-1,0],[0,0,1,0,1,1,1],           [1,0,0,1,0,0,0],[0,1,-1,1,0,-1,0],[1,0,0,1,0,0,0]])

G[801]=Matrix([[0, 1, 1, 0, -1, -1, -1],[1, 1, 1, 1, 0, 0, 0],[1, 1, 1, 1, 0, 0, 
0],[0, 1, 1, 0, -1, -1, -1],           [-1, 0, 0, -1, 0, 0, 0],[-1, 0, 0,-1, 
0, 0, 0],[-1, 0, 0,-1, 0, 0, 0]])

G[812]=Matrix([[0,1,1,0,0,0,0],[1,0,0,1,1,0,1],[1,0,0,1,1,0,1],[0,1,1,1,1,1,0],
[0,1,1,1,1,1,0],[0,0,0,1,1,1,0],[0,1,1,0,0,0,0]])

G[831]=Matrix([[0,-1,1,0,0,0,0],[-1,0,1,0,0,1,-1],[1,1,-1,1,1,0,1],
[0,0,1,1,1,1,0],[0,0,1,1,1,1,0],[0,1,0,1,1,1,0],[0,-1,1,0,0,0,0]])

G[832]=Matrix([[-1,-1,-1,0,0,0,-1],[-1,0,-2,-1,0,-1,-1],[-1,-2,-1,0,-1,0,-1],[0,-
1,0,0,-1,0,0],[0,0,-1,-1,-1,-1,0],[0,-1,0,0,-1,0,0],[-1,-1,-1,0,0,0,-1]])

G[846]=matrix([[0, -2, 0, 0, 1, 0, 0], 
[-2, 0, 1, 0, 0, 1, 1],
[0, 1,1/2, 1, 0, 1, 0],
[0, 0, 1, 2, 1, 2, 0],
[1, 0, 0, 1, 1/2,1/2, -1/2],
[0, 1, 1, 2, 1/2, 2, 0],
[0, 1, 0, 0, -1/2, 0, 0]])

G[863]=Matrix([[-1,1,1,0,0,0,-1],[1,-1,0,1,0,1,1],[1,0,0,1,1,0,1],
[0,1,1,1,1,0,0],[0,0,1,1,0,1,0],[0,1,0,0,1,-1,0],[-1,1,1,0,0,0,-1]])

G[873]=Matrix([[1,-1,-1,0,0,0,1],[-1,0,0,-1,0,-1,-1],[-1,0,0,-1,0,-1,-1],[0,-1,-
1,0,-1,0,0],[0,0,0,-1,1,-1,0],[0,-1,-1,0,-1,0,0],[1,-1,-1,0,0,0,1]])

G[878]=Matrix([[-1, 1, 0, 0, 1,0,-1],[1,-1,1,0,0,1,1],[0,1,-1,1,0,-1,0],[0,0,1,0,1,1,0],           [1,0,0,1,0,0,1],[0,1,-1,1,0,-1,0],[-1, 1, 0, 0, 1,0,-1]])

G[913]=Matrix([[1,1,1,1,0,0,0],[1,1,1,1,0,0,0],[1,1,2,2,1,1,1],[1,1,2,2,1,1,1],[0,0,1,1,0,0,0],[0,0,1,1,0,0,0],[0,0,1,1,0,0,0]])

G[918]=Matrix([[0,-1,-1,0,0,0,0],[-1,-1,-1,-1,-1,0,-1],[-1,-1,-1,-1,-1,0,-1],[0,-1,-1,-1,-1,1,0],[0,-1,-1,-1,-1,1,0],[0,0,0,1,1,-1,0],[0,-1,-1,0,0,0,0]])

G[924]=Matrix([[0,0,1,1,0,0,0],[0,0,1,1,0,0,0],[1,1,1,1,1,0,0],[1,1,1,2,2,1,1],
[0,0,1,2,1,1,1],[0,0,0,1,1,1,1],[0,0,0,1,1,1,1]])

G[932]=Matrix([[0,-1,-1,0,0,0,0],[-1,-2,-1,2,-1,0,-1],[-1,-1,-1,1,0,-1,-1],
[0,2,1,-1,1,-1,0],[0,-1,0,1,-1,1,0],[0,0,-1,-1,1,-1,0],[0,-1,-1,0,0,0,0]])

G[944]=Matrix([[1,1,1,0,0,0,1],[1,1,2,0,1,1,2],[1,2,1,1,0,0,1],[0,0,1,0,1,1,1],
[0,1,0,1,0,0,0],[0,1,0,1,0,0,0],[1,2,1,1,0,0,1]])

G[953]=Matrix([[-2,1,0,1,0,0,0],[1,0,1/2,0,1,1/2,1/2],
[0,1/2,1,1,1,0,0],[1,0,1,1/2,1,0,0],
[0,1,1,1,2,1,1],[0,1/2,0,0,1,1,1],[0,1/2,0,0,1,1,1]])

G[956]=Matrix([[1,1,0,1,1,1,1],[1,0,1,0,0,1,0],[0,1,0,1,0,0,0],[1,0,1,0,0,1,0],
[1,0,0,0,1,1,1],[1,1,0,1,1,1,1],[1,0,0,0,1,1,1]])

G[958]=Matrix([[0,1,1,-2,1,1,0],[1,1,1,0,0,0,1],[1,1,1,0,0,0,1],[-2,0,0,-2,0,0,-
2],[1,0,0,0,1,1,1],[1,0,0,0,1,1,1],[0,1,1,-2,1,1,0]])

G[970]=Matrix([[1,1,1,0,0,0,0],[1,2,2,1,1,0,0],[1,2,1,1,0,-1,-1],[0,1,1,1,1,0,0],[0,1,0,1,0,-1,-1],[0,0,-1,0,-1,-1,-1],[0,0,-1,0,-1,-1,-1]])

G[990]=Matrix([[1/8,1/2,0,1/2,0,0,0],[1/2,0,2,1,1,0,0],[0,2,-1,1,0,1,1],[1/2,1,1,3/2,1/2,0,0],[0,1,0,1/2,1/2,1,1],[0,0,1,0,1,1,1],[0,0,1,0,1,1,1]])

G[995]=Matrix([[1,1,0,1,1,0,0],[1,1,1,0,0,1,1],[0,1,-1,1,2,0,0],[1,0,1,0,-1,0,0],[1,0,2,-1,-2,1,1],[0,1,0,0,1,1,1],[0,1,0,0,1,1,1]])

G[996]=Matrix([[1,0,1,0,1,0,1],[0,2,1,1,0,2,0],[1,1,1,0,2,1,1],[0,1,0,0,1,1,0],[1,0,2,1,-1,0,1],[0,2,1,1,0,2,0],[1,0,1,0,1,0,1]])

G[1002]=Matrix([[1,0,1,1,0,0,1],[0,0,1,1,-1,-1,0],[1,1,1,0,0,-1,1],[1,1,0,-1,1,0,1],[0,-1,0,1,0,1,0],[0,-1,-1,0,1,2,0],[1,0,1,1,0,0,1]])

G[1005]=matrix(
[(1, 0, 0, 0, 1, 1, 1), 
(0, 1, 0, 1, 0, -1, 1), 
(0, 0, 1, 1, 1, 0, -1),
(0, 1, 1, 2, 1,-1, 0), 
(1, 0, 1,1, 2, 1, 0), 
(1,-1, 0, -1, 1, 2, 0), 
(1,1, -1, 0, 0, 0, 3)])

G[1028]=Matrix([[1,0,2,-1,1,1,0],[0,0,1,0,1,0,0],[2,1,4,-1,2,2,1],[-1,0,-1,1,0,-1,0],[1,1,2,0,1,1,1],[1,0,2,-1,1,1,0],[0,0,1,0,1,0,0]])

G[1060]=Matrix([[2,-1,-1,0,1,1,0],[-1,2,1,1,1,0,1],[-1,1,1,0,0,0,1],[0,1,0,1,1,0,0],[1,1,0,1,2,1,1],[1,0,0,0,1,1,1],[0,1,1,0,1,1,2]])

G[1075]=Matrix([[0,-1,2,1,2,0,0],[-1,0,1,0,0,1,-1],[2,1,0,1,0,0,2],[1,0,1,1,1,0,1],[2,0,0,1,1,-1,2],[0,1,0,0,-1,1,0],[0,-1,2,1,2,0,0]])


G[1077]=Matrix([[0,0,1,1,1,0,0],[0,1,0,-1,-1,0,1],[1,0,1,1,0,1,0],[1,-1,1,2,1,1,-1],[1,-1,0,1,0,1,-1],[0,0,1,1,1,0,0],[0,1,0,-1,-1,0,1]])

G[1087]=Matrix([[0,0,1,1,1,0,0],[0,1,1,0,1,1,0],[1,1,0,-1,0,1,1],[1,0,-1,-1,-1,0,1],[1,1,0,-1,0,1,1],[0,1,1,0,1,1,0],[0,0,1,1,1,0,0]])

G[1095]=matrix(
[(2,1,0,0,1,-1,-1), 
(1,2,1,0,0,-1,2), 
(0,1,1,1,0,0,2),
(0,0,1,3,1,1,1),
(1,0,0,1,1,0,-1),
(-1,-1,0,1,0,1,0),
(-1,2,2,1,-1,0,5)])

G[1099]=Matrix([[1,0,1,0,1,0,1],[0,2,1,1,0,2,0],[1,1,1,0,2,1,2],[0,1,0,0,1,1,1],[1,0,2,1,-1,0,-1],[0,2,1,1,0,2,0],[1,0,2,1,-1,0,-1]])

G[1104]=matrix(
[(-1,2,2,3,0,0,0), 
(2,-1,2,0,3,0,3), 
(2,2,-1,0,0,3,-3),
(3,0,0,-1,2,2,0), 
(0,3,0,2,-1,2,-3), 
(0,0,3,2,2,-1,3), 
(0,3,-3,0,-3,3,-6)])

G[1146]=matrix([[-1,1,0,0,-1,0,-3],[1,-2,1,0,1,-1,4],
[0,1,-2,1,1,1,-4],[0,0,1,-1,-1,0,3],[-1,1,1,-1,-2,0,0],[0,-1,1,0,0,-1,1],[-3,4,-4,3,0,1,-19]])


G[1167]=Matrix([[0,1,0,0,1,2,1],[1,-1,1,0,0,-1,-1],[0,1,-1,1,0,1,2],[0,0,1,-
1,1,1,-1],[1,0,0,1,0,0,1],[2,-1,1,1,0,-1,0],[1,-1,2,-1,1,0,-2]])

G[1205]=matrix(	
[(3/5, 1/10, -1/2,0,1,1,0), 
(1/10, 3/5, 1/2,1,0,1,0), 
(-1/2, 1/2, 1,1,-1,0,0),
(0, 1, 1,3, 1, 1,3),
(1, 0, -1,1, 3, 1,3),
(1, 1, 0,1, 1, 3,-1),
(0,0,0,3,3,-1,7)])


G[1212]=Matrix([(0, 4, -2, 3, 2, 1, 0), (4, -1, 1, 0, 0, 1, 4), (-2, 1,
-1, 1, 0, 0, -2), (3, 0, 1, -1, 1, 0, 3), (2, 0, 0, 1, 0, 1, 2), (1,
1, 0, 0, 1, 0, 1), (0, 4, -2, 3, 2, 1, 0)])

s = "\n".join(["%s&$\sage{G[%s]}$&\\raisebox{-1in}{\sageplot[width=2in]{plot(Graph(G[%s]), vertex_colors={'white': range(G[%s].nrows())}),figsize=[2,2]}}\\\\\\hline"%(i,i,i,i) for i in sorted(G.keys())])
print s
\end{sagesilent}

\begin{longtable}{|c|c|c|}
  \caption{Selected optimal matrices}\label{tab:optimal}\\
\hline
Atlas \#  & Matrix & Graph\\\hline\hline\endfirsthead
\hline
Atlas \#  & Matrix & Graph\\\hline\hline\endhead
\hline\endlastfoot

721&$\sage{G[721]}$&\raisebox{-1in}{\sageplot[width=2in]{plot(Graph(G[721]), vertex_colors={'white': range(G[721].nrows())}),figsize=[2,2]}}\\\hline
801&$\sage{G[801]}$&\raisebox{-1in}{\sageplot[width=2in]{plot(Graph(G[801]), vertex_colors={'white': range(G[801].nrows())}),figsize=[2,2]}}\\\hline
812&$\sage{G[812]}$&\raisebox{-1in}{\sageplot[width=2in]{plot(Graph(G[812]), vertex_colors={'white': range(G[812].nrows())}),figsize=[2,2]}}\\\hline
831&$\sage{G[831]}$&\raisebox{-1in}{\sageplot[width=2in]{plot(Graph(G[831]), vertex_colors={'white': range(G[831].nrows())}),figsize=[2,2]}}\\\hline
832&$\sage{G[832]}$&\raisebox{-1in}{\sageplot[width=2in]{plot(Graph(G[832]), vertex_colors={'white': range(G[832].nrows())}),figsize=[2,2]}}\\\hline
846&$\sage{G[846]}$&\raisebox{-1in}{\sageplot[width=2in]{plot(Graph(G[846]), vertex_colors={'white': range(G[846].nrows())}),figsize=[2,2]}}\\\hline
863&$\sage{G[863]}$&\raisebox{-1in}{\sageplot[width=2in]{plot(Graph(G[863]), vertex_colors={'white': range(G[863].nrows())}),figsize=[2,2]}}\\\hline
873&$\sage{G[873]}$&\raisebox{-1in}{\sageplot[width=2in]{plot(Graph(G[873]), vertex_colors={'white': range(G[873].nrows())}),figsize=[2,2]}}\\\hline
878&$\sage{G[878]}$&\raisebox{-1in}{\sageplot[width=2in]{plot(Graph(G[878]), vertex_colors={'white': range(G[878].nrows())}),figsize=[2,2]}}\\\hline
913&$\sage{G[913]}$&\raisebox{-1in}{\sageplot[width=2in]{plot(Graph(G[913]), vertex_colors={'white': range(G[913].nrows())}),figsize=[2,2]}}\\\hline
918&$\sage{G[918]}$&\raisebox{-1in}{\sageplot[width=2in]{plot(Graph(G[918]), vertex_colors={'white': range(G[918].nrows())}),figsize=[2,2]}}\\\hline
924&$\sage{G[924]}$&\raisebox{-1in}{\sageplot[width=2in]{plot(Graph(G[924]), vertex_colors={'white': range(G[924].nrows())}),figsize=[2,2]}}\\\hline
932&$\sage{G[932]}$&\raisebox{-1in}{\sageplot[width=2in]{plot(Graph(G[932]), vertex_colors={'white': range(G[932].nrows())}),figsize=[2,2]}}\\\hline
944&$\sage{G[944]}$&\raisebox{-1in}{\sageplot[width=2in]{plot(Graph(G[944]), vertex_colors={'white': range(G[944].nrows())}),figsize=[2,2]}}\\\hline
953&$\sage{G[953]}$&\raisebox{-1in}{\sageplot[width=2in]{plot(Graph(G[953]), vertex_colors={'white': range(G[953].nrows())}),figsize=[2,2]}}\\\hline
956&$\sage{G[956]}$&\raisebox{-1in}{\sageplot[width=2in]{plot(Graph(G[956]), vertex_colors={'white': range(G[956].nrows())}),figsize=[2,2]}}\\\hline
958&$\sage{G[958]}$&\raisebox{-1in}{\sageplot[width=2in]{plot(Graph(G[958]), vertex_colors={'white': range(G[958].nrows())}),figsize=[2,2]}}\\\hline
970&$\sage{G[970]}$&\raisebox{-1in}{\sageplot[width=2in]{plot(Graph(G[970]), vertex_colors={'white': range(G[970].nrows())}),figsize=[2,2]}}\\\hline
990&$\sage{G[990]}$&\raisebox{-1in}{\sageplot[width=2in]{plot(Graph(G[990]), vertex_colors={'white': range(G[990].nrows())}),figsize=[2,2]}}\\\hline
995&$\sage{G[995]}$&\raisebox{-1in}{\sageplot[width=2in]{plot(Graph(G[995]), vertex_colors={'white': range(G[995].nrows())}),figsize=[2,2]}}\\\hline
996&$\sage{G[996]}$&\raisebox{-1in}{\sageplot[width=2in]{plot(Graph(G[996]), vertex_colors={'white': range(G[996].nrows())}),figsize=[2,2]}}\\\hline
1002&$\sage{G[1002]}$&\raisebox{-1in}{\sageplot[width=2in]{plot(Graph(G[1002]), vertex_colors={'white': range(G[1002].nrows())}),figsize=[2,2]}}\\\hline
1005&$\sage{G[1005]}$&\raisebox{-1in}{\sageplot[width=2in]{plot(Graph(G[1005]), vertex_colors={'white': range(G[1005].nrows())}),figsize=[2,2]}}\\\hline
1028&$\sage{G[1028]}$&\raisebox{-1in}{\sageplot[width=2in]{plot(Graph(G[1028]), vertex_colors={'white': range(G[1028].nrows())}),figsize=[2,2]}}\\\hline
1060&$\sage{G[1060]}$&\raisebox{-1in}{\sageplot[width=2in]{plot(Graph(G[1060]), vertex_colors={'white': range(G[1060].nrows())}),figsize=[2,2]}}\\\hline
1075&$\sage{G[1075]}$&\raisebox{-1in}{\sageplot[width=2in]{plot(Graph(G[1075]), vertex_colors={'white': range(G[1075].nrows())}),figsize=[2,2]}}\\\hline
1077&$\sage{G[1077]}$&\raisebox{-1in}{\sageplot[width=2in]{plot(Graph(G[1077]), vertex_colors={'white': range(G[1077].nrows())}),figsize=[2,2]}}\\\hline
1087&$\sage{G[1087]}$&\raisebox{-1in}{\sageplot[width=2in]{plot(Graph(G[1087]), vertex_colors={'white': range(G[1087].nrows())}),figsize=[2,2]}}\\\hline
1095&$\sage{G[1095]}$&\raisebox{-1in}{\sageplot[width=2in]{plot(Graph(G[1095]), vertex_colors={'white': range(G[1095].nrows())}),figsize=[2,2]}}\\\hline
1099&$\sage{G[1099]}$&\raisebox{-1in}{\sageplot[width=2in]{plot(Graph(G[1099]), vertex_colors={'white': range(G[1099].nrows())}),figsize=[2,2]}}\\\hline
1104&$\sage{G[1104]}$&\raisebox{-1in}{\sageplot[width=2in]{plot(Graph(G[1104]), vertex_colors={'white': range(G[1104].nrows())}),figsize=[2,2]}}\\\hline
1146&$\sage{G[1146]}$&\raisebox{-1in}{\sageplot[width=2in]{plot(Graph(G[1146]), vertex_colors={'white': range(G[1146].nrows())}),figsize=[2,2]}}\\\hline
1167&$\sage{G[1167]}$&\raisebox{-1in}{\sageplot[width=2in]{plot(Graph(G[1167]), vertex_colors={'white': range(G[1167].nrows())}),figsize=[2,2]}}\\\hline
1205&$\sage{G[1205]}$&\raisebox{-1in}{\sageplot[width=2in]{plot(Graph(G[1205]), vertex_colors={'white': range(G[1205].nrows())}),figsize=[2,2]}}\\\hline
1212&$\sage{G[1212]}$&\raisebox{-1in}{\sageplot[width=2in]{plot(Graph(G[1212]), vertex_colors={'white': range(G[1212].nrows())}),figsize=[2,2]}}\\\hline
\end{longtable}

\end{document}